\documentclass[oneside,11pt]{amsart}
\usepackage{amssymb,mathrsfs,amstext, enumerate, comment}
\usepackage[noadjust]{cite}
\usepackage[plainpages=false,colorlinks,hyperindex,pdfpagemode=None,bookmarksopen,linkcolor=blue,citecolor=blue,urlcolor=blue]{hyperref}

\newtheorem{theorem}{Theorem}[section]

\newtheorem{corollary}[theorem]{Corollary}
\newtheorem{lemma}[theorem]{Lemma}

\theoremstyle{definition}

\theoremstyle{remark}

\def\la{\langle}
\def\ra{\rangle}

\def\Span{{\rm span}\,}

\def\cB{{\mathcal B}}

\def\cH{{\mathcal H}}

\def\Span{{\rm span}\,}

\def\diag{{\rm diag}\,}
\def\tr{{\rm tr}\,}

\def\qed{\hfill\vbox{\hrule width 6 pt\hbox{\vrule height 6 pt width 6
pt}}\medskip}

\def\IC{{\mathbb C}}

\def\IR{{\mathbb R}}

\begin{document}
\openup .5\jot
\title{Norm of an operator with numerical range in a sector}
\author{Chi-Kwong Li  and
Kuo-Zhong Wang}

\begin{abstract} We refine a recent result of Drury concerning the
optimal ratio between the norm and numerical radius
of a bounded linear operator $T$ with numerical range lying in a sector of a circular disk.
In particular, characterization is given to the operators attaining the optimal ratio, and
properties of such operators are explored.
\end{abstract}
\maketitle

Keywords. Norm, numerical radius.

AMS classification. 15A60, 47A12.

\section{Introduction}

Let $\cB(\cH)$ be the set of bounded linear operators acting on the Hilbert space $\cH$ of dimension larger than 1.
If $n$ has dimension $n$, $\cB(\cH)$ is identified with the set $M_n$ of $n\times n$ complex matrices.
For $T \in \cB(\cH)$, its numerical range and numerical radius are   defined by
$$W(T)  = \{ \la Tx,x\ra: x \in \cH, \|x\| = 1\} \qquad \hbox{ and } \qquad
w(T) = \sup\{|\mu|: \mu \in W(T)\},$$
respectively.
It is known that the numerical radius is a norm on $\cB(\cH)$ such that
$$w(T) \le \|T\| \le 2 w(T).$$
Let $\alpha \in [0, \pi/2]$ and
\begin{equation}\label{def-salpha}
S(\alpha) = \{ a+ib \in \IC: a, b \in \IR, {|b|} \le a \tan \alpha\}.
\end{equation}
Suppose $W(T) \subseteq S(\alpha)$.
In \cite{D} the author proved that
\begin{equation}\label{D-best}
\|T\| \le w(T) \sqrt{1+\sin^2\alpha}
\end{equation}
and showed that the equality is attainable. This improves the result in \cite[Theorem 3.1]{SKS} asserting $\|T\| \le w(T) \sqrt{1+2\sin^2\alpha}$
for any $T \in \cB(\cH)$.
Moreover, it was shown that the bound is attained by the matrix $A \in M_2$ equal to
$$
\frac{1}{1+2\sin^2\alpha} \begin{pmatrix}
\sqrt{1+\sin^2\alpha\cos^2\alpha} + i \sqrt{\sin^2\alpha+\sin^4\alpha} & 2\sin^2\alpha
\cr 0 & \sqrt{1+\sin^2\alpha\cos^2\alpha}-i \sqrt{\sin^2\alpha+\sin^4\alpha}\cr
\cr\end{pmatrix}.
$$
(Note that there is a small typo in the description of $A$ in the paper.)
The key step of the proof in \cite{D} used
a result concerning the characterization of a matrix $A\in M_n$ such that $w(A) \le 1$ in \cite{DW},
and some intricate  computational arguments after reducing the problem to the $2\times 2$ case.

\medskip
In this paper, we will give an elementary proof for the optimal inequality (\ref{D-best}) for $2\times 2$ matrices.
Moreover, we show that a matrix $T \in M_2$ with $W(T) \subseteq S(\alpha)$ attaining equality in (\ref{D-best}) if
and only if $T/\|T\|$ is unitarily similar to the matrix $A$ described in the preceding paragraph. These will be done in
Section 2.
As shown in \cite{D}, one can  deduce (\ref{D-best}) for a general operator $T$ once we have the result
for matrices in $M_2$. In Section 3, we  study and characterize those operators $T\in B(\cH)$ with $W(T) \subseteq S(\alpha)$
attaining the equality in (\ref{D-best}).

\medskip
For $\alpha \in [0, \pi/2]$, we always assume that $S(\alpha)$ is defined as in (\ref{def-salpha}) and
\begin{equation}\label{g-alpha}
\tau(\alpha) = \sup \{\|T\|/w(T): T \in \cB(\cH), T \ne 0, W(T) \subseteq S(\alpha)\}.
\end{equation}
If $\alpha = 0$ and $T \in \cB(\cH)$ satisfies $W(T) \subseteq S(\alpha)$,
then $W(T) \subseteq [0, \infty)$ is positive semidefinite,
and we always have $\|T\| = w(T)$ so that ${\tau}(0) = 1$.
We will exclude this case in our discussion.

\medskip
We shall always assume that
 $S(\alpha)$ and $\tau(\alpha)$ are defined as
 in (\ref{def-salpha}) and (\ref{g-alpha}) with $\alpha \in (0, \pi/2]$.
 It is easy to see that if we write $T = H + iG$ with $H = (T+T^*)/2$ and $G = i(T^*-T)/2$,
then $W(T) \subseteq S(\alpha)$ is equivalent to
$|\cos \alpha \la Gx, x\ra| \le \sin\alpha\la Hx, x\ra$
for all unit vectors $x \in \cH$, i.e.,
$\sin\alpha H  \pm  \cos \alpha G$ is positive semidefinite.

\section{The two-by-two  case}
\setcounter{equation}{0}
The goal of this section is to prove the following.

\begin{theorem} \label{main}
Let $\alpha \in (0, \pi/2]$ and $T \in M_2$ be nonzero with $W(T) \subseteq S(\alpha)$. Then
$$\|T\| \le w(T) \sqrt{1+\sin^2\alpha}.$$
The equality holds if and only if $T/\|T\|$ is unitarily
similar to the matrix
 \begin{equation}\label{T1}
\frac{1}{1+2s} \begin{pmatrix}
\sqrt{1+s-s^2} + i \sqrt{s+s^2} & 2s
\cr 0 & \sqrt{1+s-s^2}-i \sqrt{s+s^2}\cr
\cr\end{pmatrix} \ \hbox{ with } \  s = \sin^2\alpha.
 \end{equation}
\end{theorem}

We will use the following basic facts about numerical range and numerical radius; see \cite{WG}.
\begin{enumerate}
\item[(1)]   Let $T \in M_2$. Then $W(T)$  is an elliptical disk
with the eigenvalues $\lambda_1, \lambda_2$ as foci, and length of minor axis equal to
$\sqrt{\tr(TT^*)- |\lambda_1|^2 - |\lambda_2|^2}$.
\item[(2)]
Denote by  $\sigma(X)$ the spectrum of $X \in \cB(\cH)$.
If $T = H+iG \in \cB(\cH)$, where $H$ and $G$ are self-adjoint,
then $w(T) = \max_{\theta \in [0,2\pi]} \max \sigma(\cos\theta H + \sin\theta G)$.
 \end{enumerate}

To determine $\tau(\alpha)$ for $\alpha \in (0, \pi/2]$, we let
\begin{equation}\label{hattau}
\hat \tau(\alpha) = \sup \{ \|A\|/w(A): A \in R(\alpha)\},
\end{equation}
where $R(\alpha)$ is the set  of $A \in M_2$ such that
$\det(A) > 0$ and $A/\det(A)^{1/2}$ is unitarily similar to a matrix of the form
 \begin{equation} \label{T-form}
 \begin{pmatrix} r e^{i{\theta}} & 2c \cr 0 & e^{-i{\theta}}/r\cr\end{pmatrix} \quad  \hbox{ with } r \ge 1,
\ \theta \in [0, \alpha], \  c  = \sqrt{\sin^2\alpha-\sin^2\theta}.
\end{equation}
As we shall see in Lemma \ref{new-lem3},  $\hat\tau(\alpha)$ is attainable.
For $A \in R(\pi/2)$ it is  easy to check that
$A+A^*$ is positive semidefinite so that $W(A) \subseteq S(\pi/2)$.
For $\alpha \in (0, \pi/2)$, we have the following.

\begin{lemma} \label{new-lem1} Let $\alpha \in (0, \pi/2)$ and $A \in M_2$ be nonzero.
Then $A\in R(\alpha)$ or $A^*\in R(\alpha)$ if and only if $W(A)\subseteq S(\alpha)$
and touches each of the two boundary rays of $S(\alpha)$ at a nonzero point,
i.e., there are $r_1, r_2 > 0$ such that $r_1 e^{i\alpha}, r_2e^{-i\alpha} \in W(A)$.
\end{lemma}

\it Proof. \rm For the sufficiency, let
$A = H+iG$ with $H = H^*$ and $G = G^*$. We may replace
$A$ by $V^*AV$ for a suitable unitary $V \in M_2$ and assume that $H = \diag(a_1, a_2)$ with $a_1 \ge a_2 \ge 0$.
We claim that $a_2 > 0$.
If not, then either $W(A)$ is a non-degenerate  elliptical disk touching the imaginary axis so that
  $W(A) \not \subseteq S(\alpha)$, or $W(A)$ is a line segment with $0$ as an endpoint so that
$W(A)$ cannot touches the boundary rays of $S(\alpha)$ at two nonzero points. So, $a_1 \ge a_2 > 0$.
Let
$\tilde A = H^{-1/2} A H^{-1/2} = I_2 +i H^{-1/2}GH^{-1/2}$
also satisfies $W(\tilde A) \subseteq S(\alpha)$ and touches the two boundary rays of $S(\alpha)$.
This follows readily from the observation that
\begin{equation}\label{observe}
\{ v^*\tilde A v: v\in \IC^2, v \ne 0 \} =
\{v^* H^{1/2} A H^{1/2}v: v \in \IC^2, v \ne 0 \} = \{ \tilde v ^* A \tilde v: \tilde v \in \IC^2, \tilde v \ne 0\}.
\end{equation}
Since the normal matrix $\tilde A$ touches the boundary rays of $S(\alpha)$ at two nonzero points,
it  has eigenvalues $1\pm i\tan\alpha$.
As a result, $\det(\tilde A) = 1 + \tan^2\alpha$, and
$\det(A) = (1+\tan^2\alpha)/(a_1a_2) > 0$.
We may replace $A$ by $\gamma U^*AU$ with  $\gamma = \det(A)^{-1/2}$ and a suitable unitary $U \in M_2$ and assume that
$A$ is in the form (\ref{T-form}) with $r > 0$. We may assume that $r \ge 1$. Else, replace $A$ by $A^*$.
Now, let
$A = \tilde H+i\tilde G$ with $\tilde H = \tilde H^*$ and $\tilde G = \tilde G^*$.
Then $\tilde H^{-1/2}\tilde G \tilde H^{-1/2}$ has eigenvalues $\pm \tan \alpha$. Thus,
$$-\tan^2\alpha = \det(\tilde H^{-1/2}\tilde G \tilde H^{-1/2}) = \det(\tilde G)/\det(\tilde H) = -(c^2+\sin^2\theta)/(\cos^2\theta - c^2)$$
so that
$$\tan^2\alpha  = \frac{c^2-\cos^2\theta + 1}{\cos^2\theta-c^2} = -1 + \frac{1}{\cos^2\theta-c^2}
\quad \hbox{ i.e., } \quad
\frac{1}{\cos^2\alpha} = 1 + \tan^2\alpha = \frac{1}{\cos^2{\theta}-c^2}.$$
Hence, $c^2 = \cos^2\theta - \cos^2\alpha = \sin^2\alpha - \sin^2\theta.$

Conversely, suppose $A$ or $A^*$ belongs to $R(\alpha)$. Assume $A \in R(\alpha)$. The other case can be treated similarly.
Then  $A/\det(A)^{1/2}$ is unitarily similar to a matrix  of the form (\ref{T-form}).
We simply assume that $A$ has the form (\ref{T-form}), and let
$A = \tilde H+i\tilde G$ with $\tilde H = \tilde H^*$ and $\tilde G = \tilde G^*$.
Then $\tilde H$ has positive diagonal entries and $\det(\tilde H) = \cos^2\theta + \sin^2\theta - \sin^2\alpha  = \cos^2\alpha > 0$.
So, $\tilde H$ is positive definite. Moreover,
$\tilde H^{-1/2}\tilde G \tilde H^{-1/2}$ has the same eigenvalues as  $\tilde H^{-1}\tilde G$, which equals
$\frac{1}{\cos^2\alpha}\begin{pmatrix}  \cos\theta/r & -c \cr -c & r\cos\theta\cr \end{pmatrix}\begin{pmatrix}
r \sin\theta & -ic \cr ic & -\sin\theta/r\cr\end{pmatrix}$ having  trace 0 and
determinant $-\tan^2\alpha$.
As a result, $I_2 + i\tilde H^{-1/2}\tilde G \tilde H^{-1/2}$ has numerical range lying in $S(\alpha)$
and touches the boundary rays of $S(\alpha)$ at two nonzero points.
The matrix $A = \tilde H + i\tilde G$ has the same property
by observation (\ref{observe}). \qed

We will show that the value
$\hat \tau(\alpha)$ defined in (\ref{hattau}) equals $\sqrt{1+\sin^2\alpha}$  in the next two lemmas.

\begin{lemma} \label{new-lem2} Suppose $\alpha \in (0, \pi/2]$. If $A \in M_2$ has the form {\rm (\ref{T-form})}
with $r > 1$ and $\theta > 0$, then
there is $\hat A \in R(\alpha)$ with $\|A\|/w(A) < \|\hat A\|/w(\hat A)$.
\end{lemma}

\it Proof. \rm
Suppose $A$ has the form (\ref{T-form}) with
$r > 1$ and $\theta > 0$.
Then $\|A\|^2 + 1/\|A\|^2 = r^2 + 1/r^2 + 4(\sin^2\alpha - \sin^2\theta)$.
Consider $(\hat r, \hat \theta)$ with $\hat r \in (1, r)$ and $\hat \theta \in (0, \theta)$ such that
\begin{equation}\label{hatr}
r^2 + 1/r^2 + 4(\sin^2\alpha - \sin^2\theta) = \hat r^2 + 1/\hat r^2 + 4(\sin^2\alpha - \sin^2 \hat \theta).
\end{equation}
Let $\hat c = \sqrt{\sin^2\alpha - \sin^2\hat \theta}$ and
$$A(\hat r) = \begin{pmatrix}
\hat r e^{i \hat \theta} & 2\hat c \cr 0 & e^{-i\hat\theta}/\hat r\cr\end{pmatrix} = H(\hat r) + i  G(\hat r),$$
where $H(\hat r) = H(\hat r)^*$ and  $G(\hat r) = G(\hat r)^*$.  By (\ref{hatr}),
$\|A\|^2+1/\|A\|^2 = \|A(\hat r)\|^2+1/\|A(\hat r)\|^2$ so that $\|A\| = \|A(\hat r)\|$.
We will show that  {\it there is $A(\hat r)$ with $1 < \hat r < r$ and $0<\hat \theta<\theta$
such that $w(A(\hat r)) < w(A)$.
It will then follow that $\|A\|/w(A) < \|A(\hat r)\|/w(A(\hat r))$.
}

\medskip
To achieve our goal, let $\lambda_1(X)$ be the largest eigenvalues of a Hermitian matrix
$X \in M_n$.
Then
$$w(A( \hat r))= \max\{ \lambda_1(\Re (e^{-i\phi} A(\hat r))): \phi \in [0, \pi/2]\}$$
and
 $\lambda_1(\Re e^{-i\phi} A(\hat r))$ is the larger zero of the polynomial
\begin{eqnarray*}
f_{\hat r}(\lambda) &=& \det(\lambda I_2 - \cos  \phi H(\hat r) - \sin  \phi  G(\hat r))\\
&=& \lambda^2 - (\hat r\cos (\hat \theta- \phi) + \cos(\hat \theta +  \phi)/\hat r) \lambda
+ \cos(\hat \theta- \phi)\cos(\hat \theta + \phi) - \hat c^2.
\end{eqnarray*}
Since
$\cos(\hat \theta-\phi)\cos(\hat \theta + \phi) - \hat c^2
=  \cos^2 \hat \theta\cos^2  \phi-\sin^2 \hat \theta\sin^2 \phi-\sin^2\alpha+\sin^2\hat \theta
= \cos^2 \phi - \sin^2\alpha$,
$$f_{\hat r}(\lambda) = \lambda^2 - (\hat r\cos (\hat \theta- \phi) + \cos(\hat \theta +  \phi)/\hat r) \lambda
+ \cos^2\alpha -\sin^2\phi.$$
Let $\hat \phi \in (0, \pi/2)$ be such that
$w(A(\hat r))$ equal to the larger zero $\lambda_1(\hat r)$ of
$f_{\hat r}(\lambda)$.
That is, $w(A(\hat r))$ equals the larger zero of
$$f_{\hat r}(\lambda) =
\lambda^2 - (\hat r\cos (\hat \theta- \hat \phi) + \cos(\hat \theta +  \hat \phi)/\hat r) \lambda
+ \cos^2\alpha - \sin^2 \hat \phi.
$$

\medskip
We {\bf claim} that there is $\phi\in (\hat\phi, \pi/2)$ such that
the larger zero $\lambda_1(r)$ of
$$
f_r(\lambda)=\lambda^2 - ( r\cos ( \theta- \phi) + \cos(\theta + \phi)/r) \lambda
+ \cos^2\alpha - \sin^2 \phi
$$
satisfies  $\lambda_1(r) > \lambda_1(\hat r)$.
Since $\lambda_1(r) = \lambda_1(\Re(e^{-i\phi}A)) \le w(A)$, we  have
$w(A(\hat r)) = \lambda_1(\hat r) < w(A)$ as desired.

Note that our claim is valid  if there is  $\phi \in (\hat \phi, \pi/2)$ satisfies
$
r\cos ( \theta- \phi) + \cos(\theta + \phi)/r \ge \hat r\cos (\hat \theta-\hat \phi) + \cos(\hat \theta + \hat \phi)/\hat r,
$
i.e.,
\begin{equation}\label{star}
(r+1/r)\cos\theta\cos\phi+(r-1/r)\sin\theta\sin\phi\ge (\hat r+1/\hat r)\cos\hat \theta\cos\hat \phi+(\hat r-1/\hat r)\sin\hat \theta\sin\hat \phi
\end{equation}
for the following reason.
If $\phi \in (\hat\phi, \pi/2)$ satisfies (\ref{star}),
then $f_{\hat r}(\lambda) = f_{r}(\lambda) + \xi_1 \lambda + \xi_2$ with $\xi_1\ge 0$ and $ \xi_2> 0$.
Evidently, $\lambda_0 = \lambda_1(\hat r) > 0$ and
$f_r(\lambda_0) =  f_{\hat r}(\lambda_0) - \xi_1 \lambda_0 - \xi_2 < 0$. Thus, $\lambda_0 < \lambda_1(r)$.

By (\ref{hatr}), our choice of $(\hat r, \hat \theta)$ always gives
$
r^2+1/r^2-4\sin^2\theta=\hat r^2+1/\hat r^2-4\sin^2\hat\theta
$
so that
\begin{eqnarray*}
&&(\hat r+1/\hat r)^2\cos^2\hat \theta-(r+1/r)^2\cos^2\theta\\
&=&\hat r^2+1/\hat r^2-(r^2+1/r^2)+(r+1/r)^2\sin^2\theta-(\hat r+1/\hat r)^2\sin^2\hat\theta\\
&=&4\sin^2\hat\theta-4\sin^2\theta+(r+1/r)^2\sin^2\theta-(\hat r+1/\hat r)^2\sin^2\hat\theta\\
&=&(r-1/r)^2\sin^2\theta-(\hat r-1/\hat r)^2\sin^2\hat\theta > 0
\end{eqnarray*}
The last inequality above holds because $1<\hat r<r$ and $0<\hat\theta<\theta<\pi/2$.
Thus we have
\begin{equation} \label{(2)}
(\hat r+1/\hat r)\cos\hat \theta > (r+1/r)\cos\theta
\end{equation}
and
\begin{equation} \label{(1)}
(\hat r+1/\hat r)^2\cos^2\hat \theta+(\hat r-1/\hat r)^2\sin^2\hat\theta=( r+1/ r)^2\cos^2 \theta+( r-1/ r)^2\sin^2\theta.
\end{equation}

Now, we can prove (\ref{star}).
Let $a=(r+1/r)\cos\theta,b=(r-1/r)\sin\theta,\hat a=(\hat r+1/\hat r)\cos\hat\theta,$ and $\hat b=(\hat r-1/\hat r)\sin\hat \theta$.
By (\ref{(1)}), $a^2+b^2=\hat a^2+\hat b^2\equiv s$. There are $\eta,\hat \eta\in (0,\pi/2)$ such that
$\cos\eta=a/\sqrt{s}$ and $\cos\hat \eta=\hat a/\sqrt{s}$. Hence $\sin\eta=b/\sqrt{s}$ and $\sin\hat\eta=\hat b/\sqrt{s}$.
By (\ref{(2)}), we have $\cos\eta<\cos\hat\eta$, i.e., $\hat\eta<\eta$.
Now, given $\hat\phi\in (0,\pi/2)$, we consider two cases.

{\bf Case 1.} $\hat\phi<\eta$: In the case, we choose $\phi=\eta$. Then
\begin{eqnarray*}
\hat a\cos\hat\phi+\hat b\sin\hat\phi &\le& \sqrt{\hat{a}^2+\hat{b}^2}=\sqrt{a^2+b^2}=\sqrt{s}=a^2/\sqrt{s}+b^2/\sqrt{s}
=a\cos\phi+b\sin\phi.
\end{eqnarray*}

{\bf Case 2.} $\hat\phi\ge \eta$:
Let $g(t)=\cos\eta\cos t+\sin\eta \sin t-\cos\hat\eta\cos t-\sin\hat\eta\sin t$ for $\eta\le t<\pi/2$.
Then $g(\eta)=1-\cos\hat\eta\cos\eta-\sin\hat\eta\sin\eta>0$ and
\begin{eqnarray*}
g'(t)&=&-\cos\eta\sin t+\sin\eta\cos t+\cos\hat\eta\sin t-\sin\hat\eta\cos t\\
&=&\sin t(\cos\hat\eta-\cos\eta)+\cos t(\sin\eta-\sin\hat\eta)\\
&>&0  \qquad (\mbox{ since }\cos\hat\eta>\cos\eta).
\end{eqnarray*}
Hence $g$ is increasing in $[\eta,\pi/2)$ and then $g(t)>0$ for $t\in [\eta,\pi/2)$.
Thus, $$
g(\hat \phi)=\frac{1}{\sqrt{s}}(a\cos\hat \phi+b\sin\hat\phi-\hat a\cos\hat\phi-\hat b\sin\hat\phi)>0.$$
By continuity, there exists $\phi\in (\hat\phi,\pi/2)$ such that
$
a\cos\phi+b\sin\phi-\hat a\cos\hat\phi-\hat b\sin\hat\phi>0.
$
\qed

\begin{lemma} \label{new-lem3} Suppose $\alpha \in (0, \pi/2]$ and $A \in R(\alpha)$.
Then $\|A\| \le w(A)\sqrt{1+\sin^2\alpha}$. The equality holds if and only if $A/\|A\|$ is unitarily similar to the  matrix
in {\rm (\ref{T1})}.
\end{lemma}

\it Proof. \rm Let $A\in M_2$ satisfy the hypothesis. We may assume that $A$ is in the form (\ref{T-form}).
Then the product of the two singular values of $A$ equals $|\det(A)| = 1$.
So, $A$ has singular values $\|A\|$ and $1/\|A\|$. Thus,
$\|A\|^2 + 1/\|A\|^2 = \tr(A^*A) = r^2 + 1/r^2 + 4 c^2$ and
$\|A\|\pm 1/\|A\| = \sqrt{(r\pm 1/r)^2 + 4 c^2}$. Thus,
\begin{equation}\label{norm}
\|A\| = \frac{1}{2}[\sqrt{(r +  1/r)^2 + 4 c^2} + \sqrt{(r - 1/r)^2 + 4 c^2}].
\end{equation}

We will show that $\|A\| \le w(A) \sqrt{1+\sin^2\alpha}$.
The equality holds if and only if $A/\|A\|$ is unitarily similar to the matrix in (\ref{T1}).
By Lemma \ref{new-lem2}, we can focus on $A \in R(\alpha)$ of the form (\ref{T-form}) with $\theta = 0$ and $r = 1$.

{\bf Case 1} Suppose $\theta = 0$.
Then $A = \begin {pmatrix} r & 2\sin \alpha\cr 0 & 1/r\end{pmatrix}$ so that
$W(A)$ is an elliptical disk with foci $1/r, r$ and minor axis of length $2\sin \alpha$.
Hence, the major axis is the real line segment with endpoints
$$\frac{1}{2} [(r+1/r) \pm \sqrt{(r-1/r)^2 + 4 \sin^2\alpha} ]
\ \ \hbox{ so that } \ \ w(A) = \frac{1}{2}[(r+1/r) + \sqrt{(r-1/r)^2 + 4 \sin^2\alpha)}].$$
By (\ref{norm}) with $c = \sqrt {\sin^2\alpha - \sin^2\theta} = \sin\alpha$,
and the fact that  ${1 \ge  4/(r+1/r)^2}$,
$$w(A) \sqrt{1+\sin^2\alpha} - \|A\| =
\frac{1}{2}  \{(r+1/r)[\sqrt{1+\sin^2\alpha} - \sqrt{1+ 4{\sin^2\alpha}/(r+1/r)^2}]
$$
$$\ \hskip 1.3in  + (\sqrt{1+\sin^2\alpha}-1)\sqrt{(r-1/r)^2 + 4\sin^2\alpha}\}>0.$$
Thus, $\sqrt{1+\sin^2\alpha} > \|A\|/w(A)$.

\rm
{\bf Case 2} Suppose $\theta > 0$ and $r = 1$ so that
$A = \begin{pmatrix} e^{i\theta} & 2c \cr 0 & e^{-i\theta}\cr\end{pmatrix}$. Then
the elliptical disk $W(A)$ has foci  $e^{i\theta}$ and $e^{-i\theta}$,
semi-minor axis of length $c = \sqrt{\sin^2\alpha-\sin^2\theta}$.
The major axis has endpoints of the form $\cos\theta \pm ib$, where $\pm b$ are the eigenvalues of
$i(A^*-A)/2$ so that $b = \sin \alpha$.
The boundary of $W(A)$ is the set
$\{ (\cos \theta + c \cos\phi) + i \sin \alpha \sin\phi: \phi \in [0, 2\pi)\}$.
To determine $w(A)$, we can focus on $x = \cos \phi \in [0,1]$  and consider
$$f(x) =
(\cos \theta + cx)^2+ \sin^2\alpha(1-x^2)
\quad \hbox{ with } \quad
f'(x) = 2c(\cos\theta + cx) - 2x\sin^2 \alpha.$$
Note that $\sin\theta \ne 0$. Then  $f'(x) = 0$ implies
$$x =  \frac{c \cos\theta}{\sin^2\alpha - c^2} = \frac{c\cos\theta}{\sin^2\theta} \qquad \hbox{ as } \qquad
c^2 = \sin^2\alpha - \sin^2\theta.$$
We will consider two subcases depending according to
$c\cos\theta/\sin^2\theta > 1$ or
$c\cos\theta/\sin^2\theta \le 1$.
Let $s = \sin^2\alpha$. Then $\sin^2\theta = s-c^2$ and $\cos^2\theta = 1-s + c^2$.
Since $c\ge 0$ and $\sin^2\theta > 0$ the two conditions above  reduce to
$$c^2(1-s +c^2) > (s-c^2)^2 > 0 \quad \hbox{ or }
\quad c^2(1-s +c^2) \le (c^2-s)^2.$$
Equivalently, $s/\sqrt{1+s} <  c < \sqrt{s}$ or
$0\le c\le {s}/\sqrt{1+s}$. So, we analyze the following two cases.

\medskip

{\bf Case 2.a.}  Assume $c\cos\theta > \sin^2\theta > 0$, i.e., $s/\sqrt{1+s} <  c < \sqrt{s}$.  Then
\begin{eqnarray*}
f'(x) &=&2c(\cos\theta+cx)-2xs
> 2(\sin^2\theta+c^2x-xs)\\
&=& 2(\sin^2\theta-x\sin^2\theta)
\ge 0.  \qquad \qquad \qquad (\hbox{since } x\in [0,1]).
\end{eqnarray*}
Hence $f$ is increasing on $[0,1]$ so that $f$ attains it maximum at $x=1$ with $f(1)=(\cos\theta +c)^2$.
By (\ref{norm}), we have $\|A\| = \sqrt{1+c^2} + c$ so that
$$g(c) = \frac{\|A\|}{w(A)} = \frac{\|A\|}{ \sqrt{f(1)}} = \frac{\sqrt{1+c^2} + c}{\sqrt{1-\sin^2\theta} + c}
= \frac{\sqrt{1+c^2} + c}{\sqrt{\cos^2\alpha + c^2} + c}
= \frac{ \sqrt {1/c^2 + 1} + 1 } {\sqrt{\cos^2\alpha/c^2 + 1} + 1}.$$
Set $y = 1/c^2$ with $y \in (1/s, (1+s)/s^2)$.
Then
\begin{eqnarray*}
 (\ln g(y))'&=& \frac{1}{2}\left(\frac{1}{(\sqrt{y+1}+1)\sqrt{y+1}}-\frac{\cos^2\alpha}{(\sqrt{y\cos^2\alpha+1}+1)\sqrt{y\cos^2\alpha+1}}\right)\\
 &=&\frac{1}{2}\frac{y\cos^2\alpha+1+\sqrt{y\cos^2\alpha+1}-(y+1)\cos^2\alpha-\sqrt{y+1}\cos^2\alpha}
 {(\sqrt{y+1}+1)\sqrt{y+1}(\sqrt{y\cos^2\alpha+1}+1)\sqrt{y\cos^2\alpha+1}}\\
 &=&\frac{1}{2}\frac{1-\cos^2\alpha+\sqrt{y\cos^2\alpha+1}-\cos^2\alpha\sqrt{y+1}}{(\sqrt{y+1}+1)
 \sqrt{y+1}(\sqrt{y\cos^2\alpha+1}+1)\sqrt{y\cos^2\alpha+1}}>0.
 \end{eqnarray*}
Thus, $g(y)$ is an increasing function in $y$ so that
$g(y) < g(y_0)$ for $y \in (1/s, (1+s)/s^2)$
if $y_0 = (1+s)/s^2$, where $s \in (0,1]$.
We have
 $$g(y_0)
=\frac{\sqrt{1+s+s^2} +  s} { 1+s}
 < \sqrt{1+s}$$
 because
 $$ (1+s)^3-(\sqrt{1+s+s^2}+s)^2
 =s(2+s+s^2-2\sqrt{1+s+s^2})>0.$$
 Thus, $g(c) = \|A\|/w(A) < \sqrt{1+s}$ in this case.

\medskip
{\bf Case 2.b} Assume $c\cos\theta \le \sin^2\theta$, i.e.,
$0\le c\le s/\sqrt{1+s}$.
Since $f''(x) = 2 (c^2-  s) =  -2\sin^2\theta< 0$ on $[0,1]$, we see that
$f(x)$ attains its maximum at $x = c\cos\theta/\sin^2\theta$ with
\begin{eqnarray*}
f(x) &=& (\cos\theta+cx)^2+s(1-x^2)
 =
s\left(s \frac{\cos^2 \theta}{\sin^4 \theta} +  \left(1- \frac{c^2\cos^2\theta}{\sin^4\theta}\right)\right)\\
& = &
s\left(
s \frac{\cos^2\theta}{\sin^4 \theta} +  1- \frac{(s - \sin^2\theta)\cos^2\theta}
{\sin^4\theta}\right)
= s \left( 1+ \frac{\cos^2\theta}{\sin^2\theta} \right) = \frac{s}{\sin^2\theta}.
\end{eqnarray*}
We will show that $g(c) = \|A\|/w(A)$ attains its maximum at the
unique value $c =\sqrt{s-\sin^2\theta}$.
To this end, note that by (\ref{norm}) we have
$$g(c) = \frac{\|A\|}{w(A)} = \frac{\|A\|}{\sqrt{f(x)}}
= \frac{ (\sqrt{1+c^2} + c) \sin\theta} {\sqrt{s}}
= \frac{ (\sqrt{1+c^2} + c) \sqrt{s-c^2}} {\sqrt{s}}.$$
Then
$$g'(c) = \frac{1}{\sqrt{s}}\left\{
\left( \frac{c}{\sqrt{1+c^2} } + 1 \right)  \sqrt{s -c^2}
- (\sqrt{1+c^2} + c)
\frac{c}{\sqrt{s - c^2}}
\right\}$$
$$ =\frac{\sqrt{1+c^2} + c }{\sqrt{s}  \sqrt{s -c^2}{\sqrt{1+c^2}}}\left((s - c^2) -  c\sqrt{1+c^2}\right) = 0
$$
if and only if
$$(s - c^2)^2 =  c^2(1+c^2), \qquad \hbox{ i.e., } \quad
c = c_0 = \frac{s}{\sqrt{1+2s}} \in \left[0,\frac{{s}}{\sqrt{1+s}}\right].$$
Clearly, if
$0 \le c < c_0$, then  $g'(c) > 0$; if $c_0 < c \le {{s}}/{\sqrt{1+s}}$, then $g'(c) < 0$.
Then
$$g(c_0) =
\frac{1}{\sqrt s}
\left( \sqrt{ 1+\frac{s^2}{1+2s} } + \frac{s}{\sqrt{1+2s}} \right)
\sqrt{s-\frac{s^2}{1+2s} }
= \frac{1}{\sqrt s} \frac{1+2s}{\sqrt{1+2s}}\frac{\sqrt{s+s^2}}{\sqrt{1+2s}}
 = \sqrt{1+s} .$$
is the maximum. Thus, we have shown that $\hat\tau(\alpha) = \|A\|/w(A) = \sqrt{1+s}$ in this case.

Using the fact that $c^2 = s - \sin^2 \theta =  s^2/(1+2s)$, we see that
$$
\sin^2\theta = \frac{s + s^2}{1+2s}, \  \ \cos^2 \theta = \frac{1+s-{s^2}}{1+2s}, \ \hbox{ and } \
\|A\| = \sqrt{1+c^2} + c =
    \frac{1+s}{\sqrt{1+2s}}
 +    \frac{s}{\sqrt{1+2s}}
 = \sqrt{1+2s}
$$
 so that
$T/\|T\|$ is unitarily similar to $A/\|A\|$, which has the asserted form in
the theorem.  From our derivation, one sees that $A$ is uniquely determined by the function $f(x)$
and $g(c)$.
\qed

We are now ready to present the following.

\medskip\noindent
\bf Proof of Theorem \ref{main} \rm
First, suppose $\alpha \in (0, \pi/2)$.
By Lemma \ref{new-lem3},
$\hat \tau(\alpha) = \sqrt{1+\sin^2\alpha}$;
a matrix $A \in R(\alpha)$ satisfies $\|A\|/w(A) = \hat \tau(\alpha)$ if and only if $A/\|A\|$ is unitarily similar to the matrix
in (\ref{T1}).

Note that $\hat \tau (\alpha) \le \tau(\alpha)$.
Suppose $A \in M_2$ satisfies $W(A) \subseteq S(\alpha)$ and $\|A\|/w(A) = \tau(\alpha)$.
Clearly, $A$ cannot be normal; else $\|A\|/w(A) = 1$. So, $W(A)$ is a nondegenerate elliptical disk.
We prove that $W(A)$ touches the boundary rays at two nonzero points so that  $A$ or $A^*$ belongs to
$R(\alpha)$ by Lemma \ref{new-lem1}.
If not, then there is $\phi \in (-\alpha, \alpha)$ such that $e^{i\phi}W(A) = W(e^{i\phi}(A)) \subseteq S(\beta)$ and touches the
boundary rays of $S(\beta)$ for some $\beta < \alpha$.
But then
$$\tau(\alpha) = \|A\|/w(A) = \|e^{i\phi}A\|/w(e^{i\phi}A)
\le \hat \tau (\beta) = \sqrt{1+\sin^2\beta} < \sqrt{1+\sin^2\alpha} = \hat \tau (\alpha) \le \tau(\alpha),
$$
which is a contradiction.
So, $A$ or $A^*$ belongs to $R(\alpha)$ and $\|A\|/w(A)= \|A^*\|/w(A^*) \le \hat \tau(\alpha)$.
It follows that $\tau(\alpha) = \hat \tau(\alpha)$; $\|A\|/w(A) = \tau(\alpha)$ if and only if
$A/\|A\|$ is unitarily similar to the matrix in (\ref{T1}) by Lemma \ref{new-lem3}.

Now, suppose $\alpha = \pi/2$.  Using a limiting argument with  $\hat \alpha \rightarrow \pi/2$ with  $\hat \alpha \in (0, \pi/2)$, we see that
if $A \in M_2$ such that $A/\|A\|$ is unitarily similar to the matrix in (\ref{T1}) with $\alpha = \pi/2$, then
$A/w(A) = \sqrt 2$. Hence, $\tau(\pi/2) \ge \sqrt 2$.
Now, suppose $A \in M_2$ satisfies $W(A) \subseteq S(\pi/2)$ and $\|A\|/w(A)= \tau(\pi/2)$.
Then $A$ is not normal.
Let $A = H + iG$ with $H = H^*$ and $G = G^*$. We may replace $A$ by $V^*AV$ and assume that
$H = \diag(a_1, a_2)$ with $a_1 \ge a_2 \ge 0$. If $a_2 > 0$, then there is $\phi \in (-\pi/2, \pi/2)$ such that
$e^{i\phi} W(A) = W(e^{i\phi}A) \subseteq S(\beta)$ with $\beta < \pi/2$ and $W(e^{i\phi}A)$
touches the two boundary rays of
$S(\beta)$ at two nonzero points. But then
$$\tau(\pi/2) = \|A\|/w(A) = \|e^{i\phi}A\|/w(e^{i\phi}A)
\le \hat \tau (\beta) = \sqrt{1+\sin^2\beta} < \sqrt 2 \le \tau(\pi/2),
$$
which is a contradiction.
Thus, $a_2 = 0$. Clearly, $a_1 > 0$. Otherwise, $A = iG$ is normal and $\|A\|/w(A) = 1$.
Let $G = (g_{ij}) \in M_2$.
Then the elliptical disk $W(A)$ touches the imaginary axis at $i g_{22}$.
If $g_{22} \ne 0$, then there is $\phi \in (-\pi/2, \pi/2)$
such that
$e^{i\phi} W(A) = W(e^{i\phi}A) \subseteq S(\beta)$ with $\beta < \pi/2$ and
$W(e^{i\phi}A)$ touches the two boundary rays of $S(\beta)$ at two nonzero points.
Then $\|A\|/w(A) \le \sqrt{1+\sin^2\beta}$, which is a contradiction.
So, $g_{22} = 0$. Note that $g_{12} = -\bar g_{21}$ is not zero.
Else, $A$ is normal and $\|A\|/w(A) = 1$. So,  $\det(A) = |g_{12}|^2> 0$.
By a unitary similarity and replacing $A$ by $A^*$ if necessary, we may assume that $A/\det(A)^{1/2}$ has the form ({\ref{T-form}).
By  Lemma \ref{new-lem3}, we get the conclusion.
\qed

\section{General operators}
\setcounter{equation}{0}

As shown in \cite{D}, one can deduce (\ref{g-alpha}) for $T \in B(\cH)$ based on the result in $M_2$.
In fact, for a general operator $T \in \cB(\cH)$, consider
two unit vectors $x, y \in \cH$ and the compression $A$ of $T$ onto the linear span of $\{x,y\}$.
We have $W(A) \subseteq W(T) \subseteq S(\alpha)$. By the result in $M_2$,
$|\langle Tx,y\rangle| = |\langle Ax,y\rangle| \le w(A)\sqrt{1+\sin^2\alpha}
\le w(T) \sqrt{1+\sin^2\alpha}$.
Since
$$\|T\| = \sup \{ |\la Tx, y\ra|: x, y \in \cH, \|x\| = \|y\| = 1\},$$
we see that  $\|T\| \le w(T) \sqrt{1+\sin^2\alpha}$.   Moreover,
by the result for $M_2$, there is  $T_1 \in M_2$ such that $\|T_1\|/w(T_1)
= \tau(\alpha)$. One see that $T = T_1 \oplus 0$ will satisfy
$$W(T) \subseteq S(\alpha) \quad \hbox{ and } \quad \|T\|/w(T) = \|T_1\|/w(T_1) = \sqrt{1+\sin^2\alpha}.$$
Hence, $\tau(\alpha) = \sqrt{1+\sin^2\alpha}$.

\medskip
We study general operators $T \in B(\cH)$
such that $W(T) \subseteq S(\alpha)$ and $\|T\|/w(T) = \tau(\alpha)$ in the following.

\begin{theorem} \label{thm3.1}
Let $T \in B(\cH)$ be nonzero and satisfy  $W(T) \subseteq S(\alpha)$.
\begin{itemize}
\item[{\rm (a)}] Then $\|T\|/w(T) = \tau(\alpha)$
 if and only if $\tilde T = T/\|T\|$ satisfies $w(\tilde T) = 1/\tau(\alpha)$.
\item[{\rm (b)}]
If $T/\|T\|$ is unitarily similar to $T_1 \oplus T_2$ such that $T_1$ has the form in {\rm (\ref{T1})} and $w(T_2) \le 1/\tau(\alpha)$,
then $\|T\|/w(T) = \tau(\alpha)$.
\end{itemize}
\end{theorem}

\it Proof. \rm (a) For the necessity, if $\|T\|/w(T) = \tau(\alpha)$, then $\|\tilde T\|/w(\tilde T) = \tau(\alpha)$.
Since $\|\tilde T\| = 1$, we see that $w(\tilde T) = 1/\tau(\alpha)$.
The converse is clear.

(b) The assumption on $\tilde T$ ensures that $w(\tilde T) = 1/\tau(\alpha)$.  By (a), we get the conclusion.
\qed

The condition $w(\tilde T)=1/\tau(\alpha)$ in Theorem \ref{thm3.1} (a) can be checked by showing that
$$\sup \{ \max(\sigma(\Re(e^{i\theta} \tilde T)): \theta \in [-\alpha, \alpha]\} = 1/\tau(\alpha).$$
It is known that if $T \in M_n$ then $\|T\| \le 2w(T)$; the equality holds if and only if
$T/\|T\|$ is unitarily similar to $T_1 \oplus T_2$ where $T_1 = \begin{pmatrix} 0 & 1 \cr 0 & 0 \cr\end{pmatrix}$ with $w(T_2) \le 1/2$;
e.g., see \cite[Lemma 2]{JL}. We will show that  if $\alpha \in (0, \pi/2)$ and $T \in M_n$,
then $T/\|T\|$ is unitarily similar to $T_1 \oplus T_2$, where $T_1$ is the matrix in (\ref{T1}), i.e., the converse of Theorem \ref{thm3.1}
(b) holds.   Actually, we will show in Theorem \ref{3.1} that the same conclusion holds also for an infinite dimensional operator $T$
if there is a unit vector $x \in \cH$ such that $\|T(x)\| = \|T\|$.
But the result may fail if $\alpha = \pi/2$ as shown in Theorems \ref{3.3} and \ref{3.8}.

\medskip
We will always use the result in the last section asserting that the following matrix $A$ satisfies $W(A) \subseteq S(\alpha)$ and
$\|A\|/w(A) = \tau(A)$:
\begin{equation} \label{3.0a}
A = \begin{pmatrix} e^{i\theta} & 2c \cr 0 & e^{-i\theta}\cr\end{pmatrix}\quad \hbox{ with } \det(A) = 1 \ \hbox{ and } \ \|A\| = \sqrt{1+2s},\end{equation}
where
\begin{equation}\label{3.0}
0 \le \theta \le \alpha \le \pi/2, \
s = \sin^2\alpha, \  c = \frac{s}{\sqrt{1+2s}}, \
\sin^2\theta = \frac{s+s^2}{1+2s}, \
\cos^2 \theta  = \frac{1+s-s^2}{1+2s}.
\end{equation}

\begin{lemma}\label{lemB} Let
$A$ be the matrix in {\rm (\ref{3.0a})} with  $s, c, \theta$ satisfying {\rm (\ref{3.0})}. Then
$A$ is unitarily similar to
\begin{equation}\label{matrixB}
B =  \begin{pmatrix} \cos \theta + c & \sin\alpha \cr  -\sin\alpha & \cos\theta-c\end{pmatrix},
\end{equation}
and $P^t(E_{11}-E_{22})BP = \begin{pmatrix} \|A\| & 0 \cr 0 & -1/\|A\|\cr\end{pmatrix}$,
where $P$ is the orthogonal matrix with columns
$x = \begin{pmatrix} x_1\cr x_2\cr\end{pmatrix}= \xi\begin{pmatrix}\sin \alpha\cr \sqrt{1+c^2} - \cos\theta\cr\end{pmatrix}$
and $\tilde x = \begin{pmatrix} x_2 \cr - x_1\cr\end{pmatrix}$ for a positive number $\xi$.
\end{lemma}

\it Proof. \rm
Note that the largest singular values $\|A\| = \sqrt{1+2s}$ by (\ref{3.0a}). By the fact that the product of
the singular values of $A$ equal to $|\det(A)| = 1$, we see that the smaller singular value of $A$ is
$1/\|A\|$. Using the trace, determinant, and $\tr(BB^*)$, we see that $B$ is unitarily similar to the triangular matrix $A$.
Clearly, $A$ and $B$ have the same singular values.
If $D = E_{11}-E_{22}$, then  $DB$ is Hermitian with singular values $\|A\|$ and $1/\|A\|$, $\tr(DB) = 2c$.
So, $DB$ has eigenvalues $\|A\|$ and $-1/\|A\|$.
The first row of the  rank one matrix  $(DB - \|B\|I_2)$ equals
$(\cos \theta + c -\|B\|, \sin \alpha)
= (\cos\theta - \sqrt{1+c^2}, \sin \alpha).$ Thus, $DBx = \|B\|x$ and $DB \tilde x = -\tilde x/\|B\|$. The result follows.
\qed

\begin{lemma} \label{lemB-2}
Let $\alpha \in (0, \pi/2]$ and
$T \in B(\cH)$  satisfies
$W(T) \subseteq S(\alpha)$ and $\|T\|/w(T) = \tau(\alpha)$.
Suppose $x\in \cH$ is a unit vector such that
$\|T(x)\| = \|T\|$.
Then the operator matrix $T/\|T\|$ on $V \oplus V^\perp$ with $V = \Span\{x, T(x)\}$
has the form
$\tilde T = \begin{pmatrix} T_{11} & T_{12} \cr T_{21} & T_{22}\cr\end{pmatrix}$ such that
$T_{11}= B/\|B\|$, where $B$ is the matrix in {\rm (\ref{matrixB})} with
$s, c, \theta$ satisfying {\rm (\ref{3.0})}.
Moreover, the range space of $T_{12}$ lies in $\Span\{D\tilde x\}$ and
the range space of $T_{21}^*$  lies in $\Span \{\tilde x\}$, where $D=E_{11}-E_{22}$ and $\tilde x$ are defined as in Lemma {\rm \ref{lemB}}.
\end{lemma}

\it Proof. \rm
Suppose
$T_{11}\in M_2$ is the compression of $T/\|T\|$ on the subspace containing $x$ and $T(x)$. Then
$\|T_{11}\|=1$, $w(T_{11}) \le w(T/\|T\|)$ and  $W(T_{11}) \subseteq W(T/\|T\|) \subseteq S(\alpha)$.
Hence, $\tau(\alpha)=\|T\|/w(T) \le \|T_{11}\|/w(T_{11})\le  \tau(\alpha)$, and
$T_{11}$ is unitarily similar to the matrix in (\ref{T1}).
By Lemma \ref{lemB}, we may assume that $T_{11} = B/\|B\|$.

Suppose the range space of $T_{12}$ has a nonzero vector. We may adjust $\cB_1$ of the orthonormal basis of $\cB$
and assume that the leading $3\times 3$ submatrix of $\tilde T$  has the form
$\tilde B = \begin{pmatrix} B/\|B\| & u \cr v^* & b\cr\end{pmatrix}$ for some $u, v \in \IC^2$ with $u \ne 0$ and $b \in \IC$.
Now, $\tilde B$ is a submatrix of $T/\|T\|$. So, $\|\tilde B\| \le 1$.
Let $P$ be the orthogonal matrix in  Lemma \ref{lemB} and $D = E_{11}-E_{22}$. Then
$$(PD\oplus [1])\tilde B(P \oplus [1]) =\begin{pmatrix} PDBP/\|B\| & PDu\cr v^*P & b\cr\end{pmatrix},$$
where the $(1,1)$ entry of the matrix is $1$. It follows that the $(1,3)$ (and also $(3,1)$) entry of the matrix must be zero.
Thus, $u = \mu (x_2, x_1)^t$ for some $\mu \in \IC$.
Similarly, if the range space $T_{21}^*$ has a nonzero vector, we may assume that
$v^*$ in the above argument is nonzero. Then
$(PD\oplus [1])\tilde B(P\oplus [1]$ has $(1,1)$ entry equal to 1 will imply that
$v^* = \nu(x_2,-x_1)$ for some $\nu \in \IC$.
\qed

We can now prove the following theorem showing the converse of the assertion in Theorem \ref{thm3.1} (b) also holds
if $\alpha \in (0, \pi/2)$ and $T$ is norm attaining.

\begin{theorem} \label{3.1}
Let $\alpha\in (0,\pi/2)$ and $\cH$ has dimension at least $3$. Suppose $T\in B(\cH)$ is nonzero
and attains its norm at a unit vector $x\in \cH$, i.e.,   $\|T(x)\| = \|T\|$.
Then $W(T) \subseteq S(\alpha)$ and $\|T\|/w(T)=\tau (\alpha)$ if and only if
the operator matrix $T/\|T\|$ on $V \oplus V^\perp$ with $V = \Span\{x, T(x)\}$
has the form $T_1 \oplus T_2$, where   $T_1$ is unitarily similar to $B/\|B\|$
such that  $B$ is the matrix  in {\rm (\ref{matrixB})},
$W(T_2) \subseteq S(\alpha)$
and  $w(T_2) \le 1/\tau(\alpha)$.
\end{theorem}

\bf Proof. \rm  The sufficiency part follows from Theorem \ref{thm3.1} (b).
We focus on the necessity part. We may apply Lemma \ref{lemB-2}
and assume that $T/\|T\|$ has operator matrix $\tilde T$ satisfying the conclusion
of the lemma.
We will show that $T_{12}$ and $T_{21}$ are zero operators.
Suppose $T_{21}$ is not zero. By Lemma \ref{lemB-2}
and assume that  the leading $3\times 3$ matrix of $\tilde T$ has the form
 $\tilde B = \begin{pmatrix} B/\|B\|& u/\|B\| \cr v^*/\|B\| & b/\|B\|\cr\end{pmatrix}$,
where $v = \nu(x_2,-x_1)\ne 0$ and $u = \mu(x_2,x_1)^t$.
We may further replace $\tilde T$ by $F\tilde T F^*$ by $F = \diag(1,1,\xi) \oplus I$ with $\xi\nu = |\nu|$
and assume that $v^* = d(x_2, -x_1)$ with $\nu = d > 0$.

\medskip
Let $S = \diag(\frac{1}{d_1},\frac{1}{d_2} ) = \diag(\frac{1}{\sqrt{\cos\theta + c}}, \frac{1}{\sqrt{\cos\theta -c}})$.
Then $S^t B S = I_2+\tan \alpha(E_{12}-E_{21})$.
Hence
$\|B\|(S^t \oplus [1])\tilde B(S \oplus [1]) = \tilde H + i\tilde G$
with $\tilde H = \begin{pmatrix} I_2  & H_{12} \cr H_{12}^* & b_1\cr\end{pmatrix}$
and $\tilde G = \begin{pmatrix} G_{11} & G_{12} \cr G_{12}^* & b_2 \cr\end{pmatrix}$, where {$b=b_1+ib_2$},
$$H_{12} = \frac{1}{2}S(u + v),
\quad \hbox{ and } \quad G_{12} = \frac{-i}{2}S\left(u - v\right).$$
Since $W(\tilde B) \subseteq S(\alpha)$, we see that
$y^* \tilde G y \le \tan \alpha y^*\tilde Hy$ for all vector $y$.
That is
$$\tilde H - \tilde G/\tan\alpha =  \begin{pmatrix}
1 & i & y_1\cr
-i & 1 & y_2 \cr
\bar y_1 & \bar y_2 & b_1-b_2/\tan\alpha \cr\end{pmatrix}$$
is positive semidefinite, where
$$\begin{pmatrix} y_1 \cr y_2\cr\end{pmatrix} = \begin{pmatrix} \frac{x_2}{2d_1} [(\mu + d)+i (\mu - d)/\tan\alpha]\cr
\frac{x_1}{2d_2}[(\mu - d) +  i(\mu + d)/\tan\alpha]\cr
\end{pmatrix}.$$
Let $U=1/\sqrt{2}\begin{pmatrix} -i& 1 \cr 1 & -i\cr\end{pmatrix}\oplus [1]$. Then
$$
U^*(\tilde H-\tilde G/\tan\alpha)U=\left(
                                    \begin{array}{ccc}
                                      0 & 0 & \frac{iy_1+y_2}{2} \\
                                      0 & 2 & \frac{y_1+iy_2}{2} \\
                                      \frac{-i\bar{y}_1+\bar{y}_2}{2} & \frac{\bar{y}_1-i\bar{y}_2}{2} & b_1-b_2/\tan\alpha \\
                                    \end{array}
                                  \right)
$$
is positive semidefinite. Thus, $iy_1+y_2=0$.
If $y_1=y_2= 0$, then
$$0=(\mu+d)+i(\mu-d)/\tan\alpha=(\mu-d)+i(\mu+d)/\tan\alpha,$$
Hence,  $d\tan\alpha=id$,
which is a contradiction as $d>0$. So, $-iy_1=y_2\neq 0$. As a result,
$(1,1)^t$ is a multiple of
$$\begin{pmatrix} -iy_1 \cr y_2\cr\end{pmatrix} = \frac{(\mu+d)i}{2} w_1 + \frac{\mu-d}{2} w_2
\  \hbox{ with } \
w_1 = \begin{pmatrix} {-x_2}/{d_1}\cr {x_1}/{(d_2\tan \alpha)} \cr\end{pmatrix} \ \hbox{ and } \
w_2 = \begin{pmatrix} {x_2}/{(d_1\tan \alpha)}\cr {x_1}/{d_2}\cr\end{pmatrix}.$$
Equivalently,  $\begin{pmatrix} d_1 \cr d_2 \cr\end{pmatrix} = \begin{pmatrix} \sqrt{\cos \theta + c}\cr \sqrt{\cos\theta - c}\cr\end{pmatrix}$
is a multiple of $(\mu+d)i\tilde w_1 + (\mu-d)\tilde w_2$
$$
\tilde w_1 = \begin{pmatrix} {\cos \theta-\sqrt{1+c^2}}\cr \cos\alpha \cr\end{pmatrix} \ \hbox{ and } \
\tilde w_2 = \begin{pmatrix} {(\sqrt{1+c^2}-\cos\theta)}/{\tan \alpha}\cr \sin\alpha\cr\end{pmatrix}.$$
Let $\mu = \mu_1+i\mu_2$. Comparing real part and imaginary part, we see that
$(d_1, d_2)^t$ is a multiple of
$$(\mu_1+d+i\mu_2)i\tilde w_1 + (\mu_1-d+i\mu_2) \tilde w_2 = -\mu_2\tilde w_1 + (\mu_1-d) \tilde w_2  + i((\mu_1+d)\tilde w_1+\mu_2\tilde w_2).$$
So, $(\mu_1+d)\tilde w_1+\mu_2\tilde w_2 = 0$ and $(d_1, d_2)^t$ is a multiple of
$$-\mu_2 \tilde w_1 + (\mu_1-d)\tilde w_2  =
\begin{cases} \frac{d^2-\mu_1^2-\mu_2^2}{\mu_2} \tilde w_1 & \hbox{ if }\mu_2\neq 0\cr
\quad 2\mu_1\tilde w_2 &\hbox{ if }\mu_2=0.\cr\end{cases}
$$
Equivalently,
$$\frac{\cos\theta+c}{\cos\alpha}
=\frac{\sqrt{\cos\theta +c}}{\sqrt{\cos\theta - c}} = \frac{\cos\theta - \sqrt{1+c^2}}{\cos\alpha}\hbox{ or }
\frac{\sqrt{1+c^2}-\cos\theta}{\cos\alpha},$$ which is impossible because
$$
2\cos\theta=\frac{2\sqrt{1+\sin^2\alpha\cos^2\alpha}}{\sqrt{1+2\sin^2\alpha}}>\frac{1}{\sqrt{1+2\sin^2\alpha}}=\sqrt{1+c^2}-c.
$$
Hence $d=0$, which is a contradiction. Similarly, if $T_{12}$ is nonzero, we may apply the argument above with
$v = 0$ and $u = \mu(x_2, x_1)^t$ to derive a contradiction.

Thus, $\tilde T = T_1 \oplus T_2$.  Clearly, we have $W(T_2) \subseteq S(\alpha)$.
Moreover, we have $w(T_2) \le 1/\tau(\alpha)$. Otherwise, $\|T\|/w(T) = 1/w(T_2) > \tau(\alpha)$.
\qed

\medskip
Clearly, if $\cH$ has a finite dimension, then $T\in B(\cH)$ always attains its norm at a unit vector so that
the Theorem \ref{3.1} is valid. In the infinite dimensional case,
there may not exist a unit vector $x\in \cH$ such that $\|T\| = \|T(x)\|$
so that the necessity part of Theorem \ref{3.1}  may fail.
For example, we may
let $T_n\in M_2$ with $W(T_n)\subseteq S(n\alpha/(n+1))$ and $\|T_n\|/w(T_n)=1/w(T_n)=\tau(n\alpha/(n+1))$.
We let $T$ be the (infinite) direct sum of $ T_n$ for $n=1,2,\ldots$, so that $W(T)\subseteq S(\alpha)$ and
$\|T\|/w(T)=1/\lim w(T_n)=\tau(\alpha)$.
However, one can check that $T$ is not unitarily
similar to $T_1 \oplus T_2$ with  $T_1 \in M_2$ satisfying $\|T_1\|/w(T_1) = \tau(\alpha)$.
Nevertheless, we have the following.

\begin{corollary} \label{3.2}
 Suppose $T\in B(\cH)$ satisfies $W(T) \subseteq S(\alpha)$ and $\|T(x)\| \ne \|T\|$ for any unit vector $x \in \cH$.
If $\|T\|/w(T) = \tau(\alpha)$,
then there is a sequence of unitary operators $\{U_m\}$ in $B(\cH)$ such that
$$U_m^*TU_m/\|T\| = \begin{pmatrix} T_{11}(m) & T_{12}(m)\cr T_{21}(m) & T_{22}(m)\cr\end{pmatrix}\quad \hbox{ for } m = 1, 2, \dots$$
with
$T_{11}(m) \in M_2$ converging to the matrix $B/\|B\|$ in {\rm (\ref{matrixB})}, $T_{12}(m)$ and $T_{21}(m)$ converging to the zero operators.
\end{corollary}

\it Proof. \rm
Suppose $T\in B(\cH)$ satisfies $\|T\|/w(T) =\tau(\alpha)$. We may replace $T$ by $T/\|T\|$ and assume that $\|T\| = 1$.
There will be a sequence of unit vectors $\{x_m\}$ such that $\|T x_m\|$ converges to $\|T\| = 1$.
There is a sequence of unitary operators $\{U_m\}$ such that
$U_m^*TU_m = \begin{pmatrix} T_{11}(m) & T_{12}(m) \cr T_{21}(m) & T_{22}(m) \cr\end{pmatrix}$ with
$T_{11}(m) \in M_2$ satisfying $\|T_{11}(m)\|\rightarrow 1$.
Suppose $T_{21}(m)$ does not converge to 0. Then we may assume that the leading $3\times 3$ submatrix of
$U_m^*TU_m$ has the form
$\tilde B_m = \begin{pmatrix} T_{11}(m) & u_m \cr v_m^* & b_m \cr\end{pmatrix}$ so that
$\|u_m\|> \varepsilon$ for all positive integer $m$, where  $\varepsilon$ is a positive number.
Now,  $\{\tilde B_m\}$ is a bounded sequence, and has a convergent subsequence with limit of the form
$\tilde B = \begin{pmatrix} B/\|B\| & u\cr v^* & b \cr\end{pmatrix}$ such that $B/\|B\|$ has norm 1
 and $w(B/\|B\|) = 1/\tau(\alpha)$,  $\|u\| > \varepsilon$, $W(B/\|B\|)\subseteq W(\tilde B) \subseteq S(\alpha)$ and
$\|\tilde B\| /w(\tilde B) = \tau(\alpha)$. This contradicts Theorem \ref{3.1}.
Similarly, we can show that $T_{21}(m)$ converges to the zero operator.
Passing to a subsequence, we may assume that $T_{11}(m)$ converges to $T_1$, which satisfies
$W(T_1) \subseteq S(\alpha)$ and $\|T_1\|/w(T_1)=\tau(\alpha)$ by Theorem \ref{3.1}, we have the conclusion.
 \qed

Now, we turn to the case when $\alpha  = \pi/2$.
Note that $W(T) \subseteq S(\alpha)$ simply means that $T+T^*$ is positive semidefinite, and
$\tau(\pi/2) = \sqrt 2$.

\begin{theorem} \label{3.3}
Let $T$ be a $3\times 3$ nonzero matrix. Then the following conditions are equivalent:
\begin{itemize}
\item[{\rm (a)}] $W(T)\subseteq S(\pi/2)$ and $\|T\|/w(T)=\sqrt{2}$,
\item[{\rm (b)}] $T/\|T\|$ is unitarily similar to a $3\times 3$ matrix of the form
\begin{equation}
\label{3x3B}
 \left(
  \begin{array}{ccc}
    \frac{2}{3} & \frac{1}{\sqrt{3}} & d \\
    -\frac{1}{\sqrt{3}} & 0 & \sqrt{3}d \\
    d & -\sqrt{3}d & b_1+ib_2 \\
  \end{array}
\right)
\end{equation}
with
\begin{equation}\label{dbineq}
d\ge 0, \quad  b_1\ge 3d^2/2, \quad  \hbox{ and } \quad 18d^2+\sqrt{2(12d^2+b_1)^2+2b_2^2}\le  1.
\end{equation}
\end{itemize}
Moreover, if condition (b) holds with $d>0$, then $T$ is unitarily irreducible.
\end{theorem}

\it Proof. \rm   Suppose $W(T) \subseteq S(\pi/2)$ and $\|T\|/w(T) = \sqrt 2$.
By Lemma \ref{lemB-2}, we may assume that $T/\|T\|$ is unitarily similar to a matrix
of the form $\tilde T=\left(
                        \begin{array}{cc}
                          B_1 & \mu v \\
                          d\tilde v^* & b_1+ib_2 \\
                        \end{array}
                      \right)
$, where $B_1=\left(
           \begin{array}{cc}
             2/3 & 1/\sqrt{3} \\
             -1/\sqrt{3} & 0 \\
           \end{array}
         \right), v=(1,\sqrt{3})^t, \tilde v=(1,-\sqrt{3})^t, d\ge 0$, and $\mu \in \IC$.
We have $\|B_1\| = 1$ and $w(B_1) = 1/\sqrt 2$.
If  $x=(\sqrt 3,-i)^t/2$, then $\Re (e^{i\pi/4}B_1)x=1/\sqrt{2}x$.
Since $w(\tilde T)=\|\tilde T\|/\sqrt{2}=1/\sqrt{2}$, we have
$\Re (e^{i\pi/4}\tilde T)\hat x=1/\sqrt{2}\hat x$, where $\hat x=(x^t,0)^t$.
This implies that
\begin{eqnarray*}
0&=&(\bar{\mu} e^{-i\pi/4}+de^{i\pi/4})-i(\bar{\mu} e^{-i\pi/4}-de^{i\pi/4})
= de^{i\pi/4}(1+i)+\bar{\mu} e^{-i\pi/4}(1-i)\\
&=&\sqrt{2}d e^{i\pi/2}+\sqrt{2}\bar{\mu} e^{-i\pi/2}
=\sqrt{2}i(d-\bar{\mu}).
\end{eqnarray*}
Thus $\bar{\mu}=d\ge 0$ and then
$$
\Re(e^{i\theta}\tilde T)=\left(
                    \begin{array}{ccc}
                      \frac{2}{3}\cos\theta & i\frac{\sin\theta}{\sqrt{3}} & d\cos\theta \\
                     -i\frac{\sin\theta}{\sqrt{3}} & 0 & i\sqrt{3}d\sin\theta \\
                      d\cos\theta & -i\sqrt{3}d\sin\theta & b_1\cos\theta-b_2\sin\theta \\
                    \end{array}
                  \right).
$$
Since $W(\tilde T)\subseteq S(\pi/2)$, $\Re (\tilde T)=\left(
                                                                                 \begin{array}{ccc}
                                                                                   2/3 & 0 & d \\
                                                                                   0 & 0 & 0 \\
                                                                                   d & 0 & b_1 \\
                                                                                 \end{array}
                                                                               \right)
$ is positive semidefinite. Equivalently, $b_1\ge 3d^2/2$.

On the other hand, $w(\tilde T)\ge w(B_1)=1/\sqrt{2}$ and
so that $w(\tilde T)=1/\sqrt{2}$ if and only if $M_\theta =I_3/\sqrt{2}-\Re (e^{i\theta}\tilde T)$
is positive semidefinite for all $\theta\in [-\pi,\pi]$.
Clearly, the leading $1\times 1$ and $2\times 2$ principal submatrices of $M_\theta$ are positive definite for all
$\theta\in [-\pi,\pi]\setminus \{\pm \pi/4\}$.
Hence by continuity, this is also equivalent to $\det(M_\theta)\ge 0$ for all $\theta\in [-\pi,\pi]\setminus \{\pm \pi/4\}$.
Let
$\gamma = \frac{1}{\sqrt{2}}-b_1\cos\theta+b_2\sin\theta$.
By expanding cofactors along the third column of $M_\theta$, we see that $\det (M_\theta)$ equals
\begin{eqnarray*}
&&\gamma \left| \begin{matrix}
\frac{1}{\sqrt{2}}-\frac{2}{3}\cos\theta & -i\frac{\sin\theta}{\sqrt{3}} \\
i\frac{\sin\theta}{\sqrt{3}} & \frac{1}{\sqrt{2}}
\end{matrix} \right|+id\sqrt{3}\sin\theta\left| \begin{matrix}
\frac{1}{\sqrt{2}}-\frac{2}{3}\cos\theta & -i\frac{\sin\theta}{\sqrt{3}} \\
-d\cos\theta & id\sqrt{3}\sin\theta\end{matrix} \right|-d\cos\theta\left| \begin{matrix}
i\frac{\sin\theta}{\sqrt{3}} & \frac{1}{\sqrt{2}} \\
-d\cos\theta & id\sqrt{3}\sin\theta
\end{matrix} \right|\\
&=&\gamma (\frac{1}{2}-\frac{\sqrt{2}}{3}\cos\theta-\frac{\sin^2\theta}{3})
-d^2\sqrt{3}\sin^2\theta(\frac{\sqrt{3}}{\sqrt{2}}-\sqrt{3}\cos\theta)-d^2\cos\theta(-\sin^2\theta+\frac{\cos\theta}{\sqrt{2}})\\
&=&\frac{\gamma}{3}(\cos\theta-\frac{1}{\sqrt{2}})^2-d^23(1-\cos^2\theta)(\frac{1}{\sqrt{2}}-\cos\theta)
+d^2\cos\theta(\frac{1}{\sqrt{2}}-\cos\theta)(\sqrt{2}+\cos\theta)\\
&=&\frac{\gamma}{3}(\cos\theta-\frac{1}{\sqrt{2}})^2+
(\frac{1}{\sqrt{2}}-\cos\theta)4d^2(\cos^2\theta+\frac{\sqrt{2}}{4}\cos\theta-\frac{3}{4})\\
&=&\frac{\gamma}{3}(\cos\theta-\frac{1}{\sqrt{2}})^2+
(\frac{1}{\sqrt{2}}-\cos\theta)^24d^2(-\frac{3}{2\sqrt{2}}-\cos\theta)\\
&=&\frac{1}{6}(\frac{1}{\sqrt{2}}-\cos\theta)^2(\sqrt{2}-18\sqrt{2}d^2-\cos\theta(24d^2+2b_1)+2b_2\sin\theta).
\end{eqnarray*}

Let $f(\theta)=\sqrt{2}-18\sqrt{2}d^2-\cos\theta(24d^2+2b_1)+2b_2\sin\theta$. We need to determine the condition on
$(d, b_1, b_2)$ such that
$f(\theta)$ is always nonnegative. Clearly,
$f(\theta)$ is minimum when  $(\cos \theta, \sin\theta)$ is chosen such that
$\cos\theta(24d^2+2b_1)-2b_2\sin\theta = \|(24d^2+2b_1, -2b_2)\|$. Hence, we need
$$1-18d^2-\sqrt{2(12d^2+b_1)^2+2b_2^2} \ge 0 \quad \hbox{ i.e., } \quad 18d^2+\sqrt{2(12d^2+b_1)^2+2b_2^2}\le  1.$$
Conversely, we can follow the above derivation and show that if $T/\|T\|$ is unitarily similar to the matrix of the form (\ref{3x3B}),
then $T+T^*$ is positive semidefinite and $\|T\|/w(T) = \sqrt{2}$.
So, $(a)$ and $(b)$ are equivalent.

Now, suppose that $d>0$ and that $T$ is unitarily reducible. Let $T$ be unitarily similar to $T_1\oplus [r]\equiv C $, where $T_1$ is a $2\times 2$ and $r\in \mathbb{C}$.
Then $\|T_1\|=\|T\|$, otherwise $w(T)=|r|=\|T\|$, this is a contradiction. This implies that
$W(T_1)\subseteq S(\pi/2)$ and
$$
\|T\|=\|T_1\|\le \sqrt{2}w(T_1)\le \sqrt{2}w(T).
$$
Hence $\|T_1\|=\sqrt{2}w(T_1)$. By Theorem \ref{main}, $T_1$ is unitaruly similar to $1/3\left(
                                                                              \begin{array}{cc}
                                                                                1+i\sqrt{2} & 2 \\
                                                                                0 & 1-i\sqrt{2} \\
                                                                              \end{array}
                                                                            \right)
$, and $\left(
          \begin{array}{cc}
            2/3 & 1/\sqrt{3} \\
            -1/\sqrt{3} & 0 \\
          \end{array}
        \right)
$
as well. Thus, $r=b_1+ib_2$. But then $\tr (T^*T)=\tr (C^*C) <\tr (\tilde T^*\tilde T)$, which is impossible. Hence, $T$ is
unitarily irreducible.
\qed

\begin{theorem} \label{3.8}
Let $n\ge 4$ and $B$ has the form {\rm (\ref{3x3B})} with $d \in (0,  1/\sqrt{45})$.
Then there is a sufficiently small $\varepsilon>0$  such that
$$
T=\left(
    \begin{array}{cc}
      B & 0 \\
      0 & 0_{n-3} \\
    \end{array}
  \right)+ (\varepsilon + 3d^2/2 -b_1-ib_2) E_{33} + \sum_{k=1}^{n-3}\varepsilon^k (E_{(k+2)(k+3)}+E_{(k+3)(k+3)})
$$
satisfies  $W(T)\subseteq S(\pi/2)$, $\|T\|/w(T)=\sqrt{2}$, and
is unitarily irreducible.
\end{theorem}

\it Proof. \rm We may assume that $\varepsilon<1$.
Let $\hat B = \begin{pmatrix}
    \frac{2}{3} & \frac{1}{\sqrt{3}} & d \cr
    -\frac{1}{\sqrt{3}} & 0 & \sqrt{3}d \cr
    d & -\sqrt{3}d & 3d^2/2 \cr
    \end{pmatrix}$, and
$T = T_1 + T_2$ with
$$T_1 = \hat B \oplus 0_{n-3} \ \ \hbox{ and } \ \
T_2 = \varepsilon E_{33} + \sum_{k=1}^{n-3}\varepsilon^k (E_{(k+2)(k+3)}+E_{(k+3)(k+3)}).$$
Since $d \in (0,  1/\sqrt{45})$,  we can choose
$\varepsilon>0$ such that
$\hat B + \diag(0,0,\varepsilon)\equiv B_\varepsilon$
satisfies $\|B_\varepsilon\|/w(B_\varepsilon) = 1/w(B_\varepsilon) =\sqrt 2$ by Theorem \ref{3.3}.
Note that $\Re(T) = \Re(T_1) + \Re(T_2)$.
Now, $\Re(T_1) = \Re(\hat B) \oplus 0_{n-3}$
and
$\Re(T_2) = \varepsilon E_{33} +
\sum_{k=1}^{n-3}\varepsilon^kE_{(k+3)(k+3)}+\sum_{k=1}^{n-3}\frac{\varepsilon^k}{2}(E_{(k+2)(k+3)}+E_{(k+3)(k+2)})$.
Clearly, $\Re(T_1)$ is positive semidefinite;
$\Re(T_2)$ is also positive semidefinite  because  its  (real) eigenvalues
lie in the region $ \cup_{k=1}^{n-3}\{z: |z-\varepsilon^k|\le \varepsilon^k/2 + \varepsilon^{k+1}/2\} \cup \{0\}$
by the Gershgorin theorem; e.g., see \cite{varga}.  Thus, $\Re(T)$ is positive semidefinite and $W(T) \subseteq S(\pi/2)$.

Now, for any $\theta \in [0, 2\pi)$, $I_n/\sqrt 2 - \Re(e^{i\theta} T) = R_1(\theta) +  R_2(\theta)$,
where
$$R_1(\theta) = [I_3/\sqrt{2} - \Re(e^{i\theta}\hat B)-\delta E_{33}] \oplus 0_{n-3}\equiv Q(\theta)\oplus 0_{n-3}$$
and
$$
R_2(\theta) =  (\diag(0,0,\delta) \oplus I_{n-3}/\sqrt 2) - \Re(e^{i\theta} T_2),$$
where $\delta>0$ is chosen so that
the determinants of the leading $1\times 1$ and $2\times 2$ submatrices of
$Q(\theta)$ are positive and
$\det (Q(\theta))\ge 0$ for all $\theta\in [-\pi,\pi]\setminus \{\pm \pi/4\}$.
The existence of such a $\delta$ is ensured by the proof of Theorem 3.6.
Hence $Q(\theta)$ is positive semidefinite and so is $R_1(\theta)$.
We may adjust $\varepsilon$ so that $2\varepsilon <\min\{\delta, 1/\sqrt 2\}$.
Then the $(j,j)$ entry of $R_2(\theta)$ is larger than the sum of the absolute values of the off-diagonal entries
in the $j$-th row
for $j = 3, \dots, n$. Hence, $R_2(\theta)$
is positive semidefinite by Gershgorin theorem.
Thus, $w(T) \le 1/\sqrt 2$.
Since the leading $2\times 2$ submatrix of $T$ has norm 1 and numerical radius $1/\sqrt 2$, we see that $w(T) = 1/\sqrt 2$, and hence
$\|T\|/w(T) = \sqrt{2}$.

Next, we claim that $T$ is unitarily irreducible. In this case, $0\in \partial W(B_\varepsilon)$ and  $0\notin \sigma(B_\varepsilon)$.
Let $a_\varepsilon=\min\{|\lambda|:\lambda\in \sigma (B_\varepsilon)\}$.
Hence $a_\varepsilon>0$ and $$a_\varepsilon\to \min\{|\lambda|:\lambda\in \sigma (\hat B)\}>0\mbox{ as }\varepsilon\to 0^+.$$
We may assume that $\varepsilon<\min\{1,a_\varepsilon\}$.
Assume that $T$ is unitarily similar to a direct summand as $\hat T_1\oplus \hat T_2$ on $H_1\oplus H_2$.
For $k=4,\ldots,n$, let
$$
x_k=e_k+\frac{\varepsilon}{1-\varepsilon}e_{k+1}+\cdots
+\frac{\varepsilon^{(1+n-k)(n-k)/2}}{(1-\varepsilon)\cdots (1-\varepsilon^{n-k})}e_{n}.
$$
Then it is straightforward to verify that $T^*x_k=\varepsilon^{k-3}x_k$ and $\mbox{span}\{x_4,\ldots,x_n\}=\mbox{span}\{e_4,\ldots,e_n\}\equiv V_1$.
We have $\langle x_i, x_j\rangle>0$ for all $4\le i,j\le n$.
Hence $V_1$ must be in $H_1$ or $H_2$, and suppose it is in $H_1$.
We obtain that $H_2\subseteq \mbox{span}\{e_1,e_2,e_3\}\equiv V_2$
and $\{\varepsilon,\ldots,\varepsilon^{n-3}\}\subseteq \sigma (\hat T_1)$.
Assume that $x_0$ is an eigenvector of $T$ corresponding to the eigenvalue $\varepsilon$.
Since $\varepsilon<a_\varepsilon=\min\{|\lambda|:\lambda\in \sigma (B_\varepsilon)\}$,
the matrix $B_\varepsilon-\varepsilon I_3$ is invertible.
Hence the eigenvector $x_0$ can be set as $(u^t,-1,0)^t$, where $u=(B_\varepsilon-\varepsilon I_3)^{-1}(0,\varepsilon)^t$.
Thus $x_0\in H_1$ and $x_0+e_4=(u^t,0)^t\in V_2$, which implies that $H_1\cap V_2\neq \{0\}$.
As a result, $B_\varepsilon$ is unitarily reducible, which is a contradiction.
Hence $T$ is unitarily irreducible.
\qed

\begin{theorem} \label{thm:pi/2a}
Suppose $\cH$ has dimension at least $3$ and $T \in B(\cH)$ is such that $T+T^*$ is positive semidefinite.
Suppose $T$ attains its norm at a unit vector $x\in \cH$, i.e.,   $\|T(x)\| = \|T\|$.
Then $\|T\|/w(T)=\sqrt 2$ if and only if the compression of the operator $\tilde T = T/\|T\|$ onto the subspace spanned by
$x$ and $T(x)$ is unitarily similar to the matrix
$B = \begin{pmatrix} 2/3 & 1/\sqrt{3}\cr -1/\sqrt{3} & 0 \cr\end{pmatrix}$ and $w(\tilde T) \le 1/\sqrt 2$.
\end{theorem}

\it Proof. \rm To prove the necessity, assume that $T$ has the said properties. Suppose
$\tilde T = T/\|T\|$ has operator matrix $\begin{pmatrix} T_{11} & T_{12} \cr T_{21} & T_{22}\cr\end{pmatrix}$,
where $T_{11}$ is a compression of $\tilde T$ onto  $\Span \{x, T(x)\}$.
Then  $T_{11} + T_{11}^*$ is positive semidefinite, $\|T_{11}\| = 1$, and $w(T_{11}) \le w(\tilde T)$.
So, $\sqrt 2 \ge \|T_{11}\|/w(T_{11}) \ge \|\tilde T\|/w(\tilde T) = \sqrt 2$.
Hence, $T_{11}$ is unitarily similar to the matrix $B$ by Lemma \ref{lemB}.

Conversely, suppose the compression of $\tilde T$ onto $\Span\{x, Tx\}$
is unitarily similar to $B$.  Then $\|\tilde T\| = \|B\| = 1$.
Since $1/\sqrt 2 \ge w(\tilde T) \ge w(B) = 1/\sqrt 2$,  we have $\|\tilde T\|/w(\tilde T) {\color{blue}=}  \|B\|/w(B) = \sqrt 2$.
\qed

\begin{corollary}
Suppose $T\in B(\cH)$ is nonzero such that $T+T^*$ is positive semidefinite.  If $\|T\|/w(T) = \tau(\pi/2)$,
then there is a sequence of unitary operators $\{U_m\}$ in $B(\cH)$ such that
$$U_m^*TU_m/\|T\| = \begin{pmatrix} T_{11}(m) & T_{12}(m)\cr T_{21}(m) & T_{22}(m)\cr\end{pmatrix}\quad \hbox{ for } m = 1, 2, \dots$$
with
$T_{11}(m) \in M_2$ converging to the matrix $B= \begin{pmatrix} 2/3 & 1/\sqrt{3}\cr -1/\sqrt{3} & 0 \cr\end{pmatrix}$.
\end{corollary}

\medskip\noindent
{\bf \large Acknowledgment}

Li is an affiliate member of the Institute for Quantum Computing, University of Waterloo.
His research was supported by the Simons Foundation Grant  851334. The research started while Li was visiting
Taiwan in January, 2024, supported by the
Mathematics Research Promotion Center,
National Science and Technology Council, Taiwan. The research of
 Wang was partially supported by National Science and Technology Council, Taiwan, under the research grant 112-2115-M-A49-003-MY2.

\medskip\noindent
(Li) Department of Mathematics, College of William \& Mary, VA 23187, USA. ckli@math.wm.edu.

\medskip\noindent
(Wang) Department of Applied Mathematics, National Yang Ming Chiao Tung University, Hsinchu

\quad 30010, Taiwan.
kzwang@math.nctu.edu.tw
\end{document}